\documentclass{amsart}
\usepackage{epsf}
\newtheorem{lemma}{Lemma}[section]
\newtheorem{propos}[lemma]{Proposition}

\newtheorem{corol}[lemma]{Corollary}
\theoremstyle{definition}

\theoremstyle{remark}

\newcommand{\CC}{\hbox{{$\sl C$}}}

\newcommand{\F}{\mathbb{F}}
\newcommand{\Z}{\mathbb{Z}}

\newcommand{\del}{\partial}

\newcommand{\eps}{{\epsilon}}
\newcommand{\tens}{\mathop{\otimes}}

\newcommand{\id}{{\rm id}}

\newcommand{\eproof}{$\quad \diamond$\bigskip}

\newcommand{\und}[1]{{\underline {#1}}}
\newcommand{\eqn}[2]{\begin{equation}#2\label{#1}\end{equation}}

\begin{document}



\title{\rm\large BRAIDED LINE AND COUNTING FIXED POINTS OF $GL(d,\F_q)$}
\author{P.J. Cameron $+$ S. Majid}%

\address{School of Mathematical Sciences\\
Queen Mary, University of London\\ 327 Mile End Rd,  London E1
4NS, UK}

\thanks{S.M. is a Royal Society University Research Fellow}%


\maketitle

\begin{abstract}We interpret a recent formula for counting orbits of
$GL(d,\F_q)$ in terms of counting fixed points as addition in the
affine braided line. The theory of such braided groups (or Hopf
algebras in braided categories) allows us to obtain the inverse
relationship, which turns out to be the same formula but with $q$
and $q^{-1}$ interchanged (a perfect duality between counting
orbits and counting fixed points). In particular, the probability
that an element of $GL(d,\F_q)$ has no fixed points is found to be
the truncated q-exponential
\[ {}_de_q^{-{1\over q-1}}=\sum_{i=0}^d {(-1)^i\over
(q^i-1)\cdots(q-1)}.\]
\end{abstract}

\section{Introduction}

\bigskip

Let $\Sigma=\{1,\cdots,d\}$. It is known \cite{Bos} that one may
write the number $F_j$ of orbits of $G\subseteq S_d$ acting on the
set $\Sigma_j$ of distinct $j$-tuples in terms of the
probabilities $p_i$ that an element of $G$ has $i$ fixed points.
Writing generating functions
\[ F=\sum_{j=0}^d F_j {x^j\over j!},\quad p(x)=\sum_{i=0}^d p_i
x^i,\] the two are related by
\[ F(x)=p(x+1).\]
For example, when $G=S_d$, it is evident that $F_j=1$ for all $j$
and hence
\[ p_0={}_de^{-1}\equiv\sum_{i=0}^d {(-1)^i\over i!}.\]

Recently, one of us\cite{Cam:unp} has considered a similar problem
for $G\subseteq GL(d,\F_q)$ where $\F_q$ denotes the field of
prime-power order $q$. Let
\[ \Sigma_j=\{(a_1,\cdots,a_j),\quad a_i\in \F_q^d,\quad a_i\
{\rm linearly\ indept.}\},\quad L_j=\#\{{\rm Orbits\ of}\ G\ {\rm
in}\ \Sigma_j\},\] and let $P_i$ denote the probability that an
element of $G$ pointwise fixes exactly a subspace of $\F_q^d$ of
dimension $i$. Then similarly by the orbit-counting lemma one has
 \eqn{LP}{ L_j={1\over |G|}\sum_{g\in G}\#\{g-{\rm fixed\ points\
in}\ \Sigma_j\}=\sum_{i=j}^d
P_i(q^i-1)(q^i-q)\cdots(q^i-q^{j-1}).} The generating function
\eqn{Lx}{ L(x)=\sum_{j=0}^d{L_j x^j\over(q^j-1)\cdots (q-1)}} is
likewise written in terms of the $P_i$ and $q$-binomial
coefficients.  The latter are well-known (the Gauss formula) to
count the number of subspaces of a given dimension in an
$\F_q$-vector space. The $q$-factorial denominators in such
expressions are understood as $1$ for $j=0$.
\bigskip

On the other hand, for applications, we would really like the
inverse relationship $P_i$ in terms of $L_j$. In this note we
provide this inverse in a natural way. Thus:

\begin{propos} With the notations above the probabilities $P_i$
of an $i$-dimensional fixed subspace are given by
\[ P_j={q^{-{j(j-1)\over 2}}\over (q^j-1)\cdots (q-1)}\sum_{i=j}^d
{L_i (-1)^{i-j}\over (q^i-q^j)\cdots (q^{j+1}-q^j)}.\]
\end{propos}

For example, it is evident that if $G=GL(d,\F_q)$ then $L_j=1$ for
all $j$. Hence we find:

\begin{corol} For $G=GL(d,\F_q)$ the probability of no fixed
points is
\[ P_0=\sum_{i=0}^d {(-1)^i\over (q^i-1)\cdots (q-1)}\equiv
{}_de_q^{-{1\over (q-1)}}.\]
\end{corol}

Our standard notation for the (truncated) q-exponential here is
such that ${}_de^x_q\to {}_de^x$ as $q\to 1$, i.e. it is defined
with the factorial of $q$-integers $[i]_q={q^i-1\over q-1}$. In
our case $q\ge 2$ and these $(q-1)$ factors cancel but the
notation reminds us of a divergence if we were to try include the
finite group case as $q=1$.

There are of course many ways to invert a linear relationship. We
mention that one may also, for example, use the method of M\"obius
inversion with respect to a suitable function, see \cite[p.
127]{Sta:enu}. On the other hand we will obtain the above
inversion in a way that is interesting because it is precisely
analogous to the computation for the $S_d$ case based on
generating functions. The only difference is that that one should
consider the generating function dummy variable $x$ for the
purposes of addition and subtraction as a braided one, i.e. the
polynomials $k[x]$ as a Hopf algebra, but in a braided category
where $x$ has braid statistics $q$ (the affine `braided line').
(The idea here is similar to that of a superspace but instead of a
`fermionic' factor $-1$ there is a factor $q$ when $x$ is
interchanged with another independent $x$.) The systematic passage
from permutation groups to linear groups over finite fields by
means of switching on a braiding would appear to be an interesting
direction suggested by the present work.

The required methods of Hopf algebras in braided categories or
`braided groups' have been introduced some years ago\cite{Ma:bg},
precisely as the the basic structure underlying quantum groups and
$q$-analysis, \cite{Ma:introm}\cite{Ma:book}. We recall what we
need in a self-contained preliminary Section~2. Working with $x$
as a braided variable we relate, in Section~3, the signed
generating function $\bar L(x)=L(-x)$ and $P(x)=\sum_{i=0}^dP_i
x^i$ by \eqn{barLP}{ \bar L(x)=P(y-_q x)|_{y=1},} where $-_q$ on
the right is subtraction under the additive braided group
structure of the braided line. We can then immediately invert the
relationship as \eqn{PbarL}{ P(x)=\bar L(y-_{q^{-1}} x)|_{y=1},}
where $-_{q^{-1}}$ is subtraction under the opposite additive
braided group structure (it is not commutative). This gives the
above inversion formula, essentially without computation. It also
exhibits a novel duality phenomenon whereby the coefficients of
$P(x)$ are given in terms of those of $\bar L(x)$ by exactly the
same formula as vice-versa, but with $q$ and $q^{-1}$
interchanged. This is discussed in Section~4.

\section{Preliminaries: the affine braided line}

We briefly recall the definition of braided groups and the braided
line. Once these are understood, the computation itself in the
next section will be immediate.

The formal definition of a braided category is in
\cite{JoyStr:bra}. It consists of a category $\CC$, a functor
$\tens:\CC\times\CC\to \CC$ and natural transformations
$\Psi:\tens\to \tens^{\rm op}$ and $\Phi:(\ \tens\ )\tens\to
\tens(\ \tens\ )$ obey respectively hexagon and pentagon
identities. In addition there is a trivial object $\und 1$ and
associated functorial isomorophisms of identity with respect to
$\tens$. A braided group\cite{Ma:bg}\cite{Ma:introm} is first of
all a unital algebra $B$ in such a category, defined in the
obvious way except that the identity must be viewed as a morphism
$\eta:\und 1\to B$. In addition, we require a coalgebra structure
$\Delta:B\to B\tens B$ (the coproduct) and $\eps:B\to \und 1$ (the
counit) obeying axioms which are the same as for an algebra but
with arrows reversed. These $\Delta,\eps$ are required to be
algebra homomorphisms, where $B\tens B$ has the braided tensor
product algebra structure. Its product $(B\tens B)\tens (B\tens
B)\to B\tens B$ consists of $\Psi$ to interchange the middle
factors followed by products $\cdot$ in $B$, as well as $\Phi$ as
needed to order brackets. Finally, there is an antipode or
`linearized inverse' morphism $S:B\to B$ obeying \eqn{antipode}{
\cdot\circ (\id\tens S)\circ
\Delta=\cdot\circ(S\tens\id)\circ\Delta=\eta\circ\eps.}

\begin{figure}
\[{\rm (a)}\quad \epsfbox{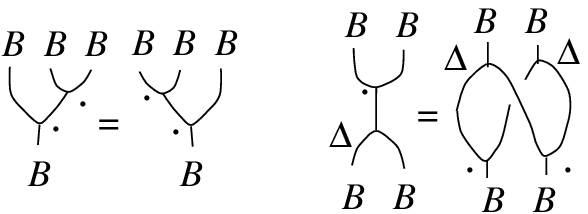},\quad {\rm (b)}\quad
\epsfbox{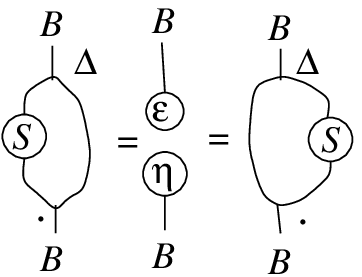}\]
\caption{Axioms of a braided group: (a) associativity of the
product (similarly for $\Delta$) and the braided homomorphism
axiom, (b) the antipode axioms}
\end{figure}

The theory of braided groups can be developed similarly to that of
quantum groups or Hopf algebras. A useful (but not necessary) tool
is to use a diagrammatic notation as follows. One omits $\Phi$ and
brackets. One denotes $\tens$ by juxtaposition, the braiding by
$\Psi=\epsfbox{braid.eps}$ and its inverse by
$\Psi^{-1}=\epsfbox{braidinv.eps}$. Algebraic operations are then
drawn as `flow charts' flowing from the top of the page downwards.
General morphisms such as $S$ are drawn as nodes. All of this is
as for any braided category. Finally, we denote the product and
coproduct of $B$ by $\Delta=\epsfbox{deltafrag.eps}$ and
$\cdot=\epsfbox{prodfrag.eps}$. The principal axioms of a braided
group are shown in this notation in Figure~1. One of the key
lemmas that we will need is:
\begin{lemma}\cite{Ma:introm} $S$ is braided-antimultiplicative
and anticomultiplicative, \[ S\circ\cdot=\cdot\circ\Psi\circ
(S\tens S),\quad \Delta\circ S=(S\tens S)\circ\Psi\circ \Delta.\]
\end{lemma}

We will also need that if $B$ is a braided group in $\CC$ with $S$
invertible then
\eqn{deltaop}{\Delta^{op}=\Psi^{-1}\circ\Delta,\quad
S^{op}=S^{-1}}
with the original algebra and counit provide
another braided group structure $B^{\rm op}$ in the category
$\CC^{\rm op}$ defined as $\CC$ but with $\Psi,\Psi^{-1}$
interchanged (i.e. with reversed braid crossings).

In our elementary application we work in the category $\CC_q$ of
$\Z$-graded vector spaces $V,W$, etc. with its usual $\tens$ and
 (trivial) $\Phi$, but with
\eqn{Psiq}{\Psi_{V,W}(v\tens w)=w\tens v\ q^{{\rm deg}(v){\rm
deg(w)}}}extended linearly. Morphisms are degree preserving maps.
We work over any field $k$. The unit object is $\und 1=k$.

The algebra of polynomials $B=k[x]$ lives in this category (all
its structure maps are morphisms) with ${\rm deg}(x)=1$. We make
it into a braided group (the `affine braided line') with
\[ \Delta x=x\tens 1+1\tens x,\quad \eps x=0,\quad Sx=-x.\]
This looks like the usual additive group Hopf algebra structure
on an affine line but $\Delta$ extends differently as a
homomorphism to the braided tensor product. The braided tensor
algebra $B\tens B$ is generated by $y\equiv x\tens 1$ and
$x\equiv 1\tens x$ with the relations $xy=qyx$ (the so-called
`quantum plane'). Then
\[ \Delta f(x)=f(y+x),\quad \eps f(x)=f(0),\]
which exhibits $\Delta$ as the operation of addition with respect
to braided variables where independent $x,y$ $q$-commute (to
emphasise this one may write $f(y+_qx)$, which is the notation
indicated in (\ref{barLP}-(\ref{PbarL})). Working directly from
the braided anti/homomorphism properties of $S$ and $\Delta$ one
also has explicitly \eqn{DeltaSx}{ \Delta
x^i=\sum_{j=0}^i\left[{i\atop j}\right]_q x^j\tens x^{i-j},\quad
S(x^i)=q^{i(i-1)\over 2}(-x)^i,} where \eqn{qbinom}{\left[{i\atop
j}\right]_q ={(q^i-1)\cdots (q^{i-j+1}-1)\over(q^j-1)\cdots
(q-1)}={[i]_q!\over [j]_q![i-j]_q!}.} Clearly, the elementary
properties of braided groups mentioned above can easily be
verified directly from these definitions and boil down to a number
of useful identities for $q$-binomial coefficients, not all of
them obvious.

The braided line can be viewed as the `correct' setting for
q-analysis in one variable including Jackson's q-derivative and
q-integration\cite{Jac:int}, q-exponentials, as well as a new
theory of q-Fourier transforms\cite{Ma:fre}\cite{KemMa:alg}. Such
methods in particular solved the long-standing problem of
extending $q$-differential calculus to the multivariable case on
affine $q$-planes (these turn out to be additive braided groups).

\section{Inverting the map from $P[x]$ to $L(x)$}

As in \cite{Cam:unp} we write (\ref{LP}) as
\[ L(x)\equiv \sum_{j=0}^d {L_j\over (q^j-1)\cdots (q-1)}
x^j=\sum_{j=0}^d\sum_{i=j}^d
P_i{(q^i-1)\cdots (q^{i-j+1}-1)\over (q^j-1)\cdots
(q-1)}q^{j(j-1)\over 2}x^j\] or, changing the order of summation
and using the q-binomial notation,
\eqn{LPbinom}{L(x)=\sum_{i=0}^dP_i\sum_{j=0}^i\left[{i\atop
j}\right]_q q^{j(j-1)\over 2}x^j .}

Now consider the algebra homomorphism $\phi:k[x]\to k$ which
`evaluates at 1', \eqn{phi}{ \phi(x^m)=1,\quad \forall m.} Using
this and the structure of the braided line recalled above, we
clearly have \begin{eqnarray} \bar
L(x)&=&\sum_{i=0}^dP_i\sum_{j=0}^i\left[{i\atop j}\right]_q
q^{j(j-1)\over 2}(-x)^j\label{barLxP}\\ &=& (\phi\tens
S)\sum_{i=0}^dP_i\sum_{j=0}^i\left[{i\atop j}\right]_q
x^{i-j}\tens x^j=(\phi\tens S)\circ\Delta P(x)\nonumber
\end{eqnarray} which is the precise meaning of (\ref{barLP}).

\begin{lemma} The linear operator $T=(\phi\tens S)\circ\Delta$
which maps $P(x)$ to $\bar L(x)$
has inverse \[ T^{-1}=(\phi\circ S\tens\id)\circ\Delta\circ
S^{-1}.\]
\end{lemma}
\begin{figure}
\[\epsfbox{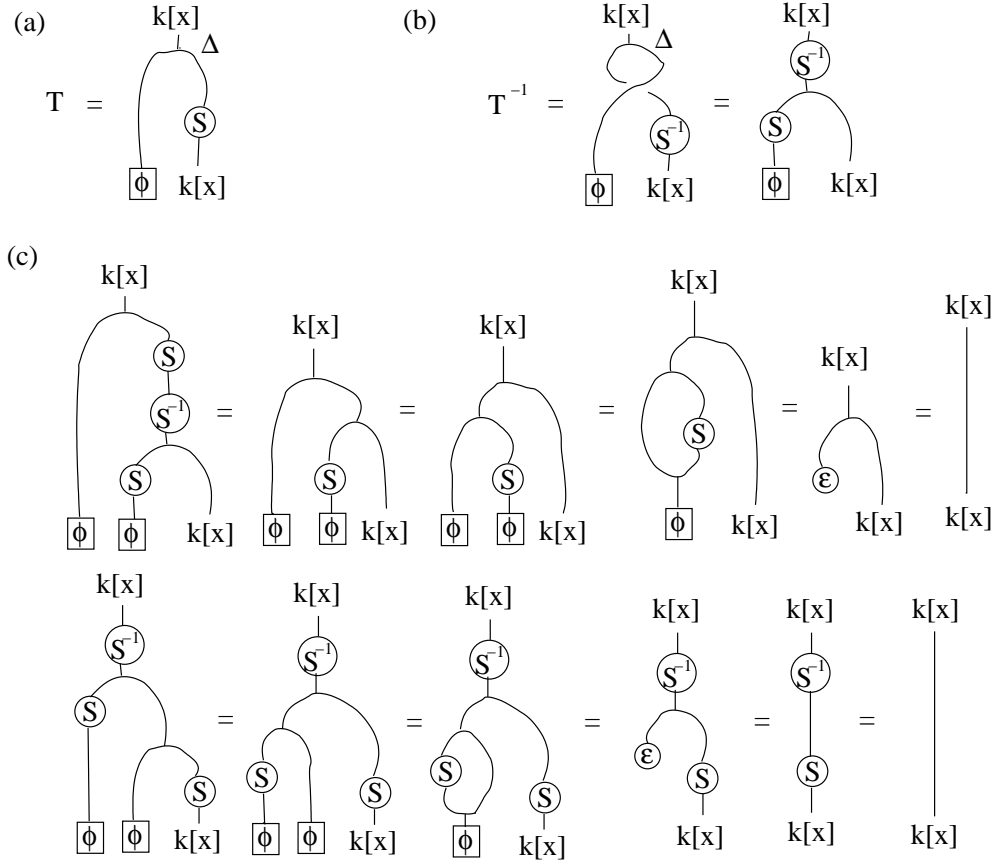}\]
\caption{(a) The operator $T$ sending $P$ to $L$ and (b) its
inverse as proven in (c)}
\end{figure}
\proof The proof is shown in Figure~2 in the diagrammatic
notation. The composition $T^{-1}\circ T$ is depicted in the first
line in Figure~2(c). Here $T^{-1}$ is linear so we can move the
output of $\phi$ in $T$ to the end of the computation as shown
(remember that juxtaposition at the bottom of the line denotes
tensor product). We cancel $S^{-1}\circ S$, then use
coassociativity of $\Delta$. Finally we use that $\phi$ is an
algebra homomorphism and cancel an `antipode loop' using the axiom
of an antipode. The resulting counit combines with $\Delta$ to the
identity (the counity axiom). Note that $\phi$ is not a morphism
in $\CC_q$ but we do not use this (we apply $\phi$ only at the end
of each diagrammatic expression). Similarly for the lower line of
Figure~2(c). Since these techniques may be unfamiliar, we also
write the proof more conventionally. Thus,
\begin{eqnarray*} T^{-1}\circ T &=& (\phi\circ
S\tens\id)\circ\Delta\circ (\phi\tens\Delta)=(\phi\tens\phi\circ
S\tens\id)\circ(\id\tens\Delta)\circ\Delta\\
&=& (\phi\tens\phi\circ S\tens\id)\circ(\Delta\tens\id)\circ\Delta
=\left(\phi\circ\cdot\circ(\id\tens
S)\circ\Delta\tens\id\right)\circ\Delta\\
& =&\phi(1)(\eps\tens\id)\circ\Delta=\id\end{eqnarray*} and
similarly for $T\circ T^{-1}$. These proofs amount to a series of
q-identities among q-binomial coefficients by writing out $\Delta,
S$ explicitly on monomials $x^m$. \eproof

It remains only to compute $T^{-1}$ explicitly. We write $\bar
L=\sum_{i=0}^d \bar L_i x^i$ (without any factorial factors). Then
\[ P(x)=T^{-1}(\bar L(x))=(\phi\circ S\tens\id)
\sum_{i=0}^d \bar L_i (-1)^iq^{-i(i-1)\over 2}\sum_{j=0}^i
\left[{i\atop j}\right]_q x^{i-j}\tens x^j\]
\[=\sum_{i=0}^d \bar L_i (-1)^iq^{-i(i-1)\over 2}\sum_{j=0}^i
(-1)^{i-j}q^{(i-j)(i-j-1)\over 2}\left[{i\atop j}\right]_q x^j
\]
giving \eqn{PbarLx}{ P(x)=\sum_{i=0}^d \bar L_i \sum_{j=0}^i
\left[{i\atop j}\right]_q q^{-j(j-1)\over 2}q^{-j(i-j)}(-x)^j.}

We can then put in the form of $\bar L_i$ in terms of $L_i$ and
reorder the summation to extract the $P_i$. Thus
\[ P(x)=\sum_{j=0}^d x^jq^{-j(j-1)\over
2}\sum_{i=j}^d L_i{(-1)^{i-j} q^{-j(i-j)}\over (q^j-1)\cdots
(q-1)(q^{i-j}-1)\cdots (q-1)}\]
\[=\sum_{j=0}^d x^j{q^{-j(j-1)\over
2} \over (q^j-1)\cdots (q-1)}\sum_{i=j}^d L_i{(-1)^{i-j} \over
(q^i-q^j)\cdots (q^{j+1}-q^j)}\] or $P_i$ as stated in
Proposition~1.1. This completes our proof of that.

\section{Duality under $q\leftrightarrow q^{-1}$}

Using one of the basic results from braided group theory, namely
the second part of Lemma~2.1, we can also write \eqn{Tinvdua}{
T^{-1}=(\phi\tens S^{-1})\circ\Psi^{-1}\circ\Delta,} which is the
other expression shown in Figure~2(b) and the precise meaning of
(\ref{PbarL}). This computes more directly as (\ref{PbarLx}), with
$q^{-j(i-j)}$ coming from the braiding $\Psi^{-1}(x^j\tens
x^{i-j})$ and $q^{-j(j-1)\over 2}$ from $S^{-1}$.

\begin{corol} The operator $T^{-1}$ is exactly the same as the operator $T$
but with $q$ replaced by $q^{-1}$.
\end{corol}
\proof We see that $T^{-1}$ has just the same form as $T$ but with
$\Delta^{\rm op}$ and $S^{\rm op}$ in place of $\Delta, S$. This
is the opposite braided group structure on $k[x]$, namely the same
form on degree 1 but developed as a braided group in the category
$\CC_q^{\rm op}=\CC_{q^{-1}}$. One may also see this explicitly by
observing that
\[\left[{i\atop j}\right]_{q^{-1}}=q^{-j(i-j)}\left[{i\atop
j}\right]_q\] so that (\ref{PbarLx}) can be written as
\eqn{PbarLx2}{P(x)=\sum_{i=0}^d \bar L_i \sum_{j=0}^i
(-1)^j\left[{i\atop j}\right]_{q^{-1}}  q^{-j(j-1)\over 2} (-x)^j}
which has the same form as our original (\ref{barLxP}). \eproof

Finally, we mention another way to obtain the above duality using
braided group theory, this time in the form of `q-analysis'.
Namely, the displacement of a function by $1$ can be written as a
q-exponential of q-derivatives (the braided Taylor's
theorem\cite{Ma:fre}). In the present one-variable case it is
easily verified on monomials and takes the form
\[ (\phi\tens\id)\Delta =e_q^{\del_q}
=\sum_{i=0}^\infty {(\del_q)^i\over [i]_q!}\]
acting on polynomials (the sum is only finite).  Here
\[ \del_qf(x)={f(qx)-f(x)\over x(q-1)}\]
is Jackson's $q$-derivative\cite{Jac:int}. (Also in our case one
can cancel the $(q-1)$ factors here and in $[i]_q$! if one
prefers.) Moreover,
\[ e_q^Ae_{q^{-1}}^{-A}=1\]
for any nilpotent operator $A$ as one may see by
q-differentiating\cite[Lem. 3.2.2]{Ma:book}. Hence
\[ T=S\circ e_q^{\del_q}=e_q^{-\del_{q^{-1}}}\circ S.\]
The second equality here is easily verified on monomials. Then
\[ T^{-1}=e_{q^{-1}}^{-\del_q}\circ S^{-1}=S^{-1}\circ
e_{q^{-1}}^{\del_{q^{-1}}},\]
which has the same form with $q$ and $q^{-1}$ interchanged.

\baselineskip 14pt

\begin{thebibliography}{10}

\bibitem{Bos}
N. Boston, W. Dabrowski, T. Foguel, P.J. Gies, J.~Leavitt,
D.T.~Ose and D.A. Jackson.
\newblock The proportion of fixed-point-free elements of a
transitive permutation group.
\newblock {\em Commun. Algebra} 21:3259--3275, 1993.

\bibitem{Cam:unp}
P.J.~Cameron.
\newblock Seminar, November 2000 (unpublished).

\bibitem{Sta:enu}
R.P.~Stanley.
\newblock {\em Enumerative Combinatorics, Vol 1}.
\newblock Cambridge University Press, 1997.

\bibitem{Ma:bg}
S.~Majid.
\newblock Braided groups.
\newblock {\em J. Pure and Applied Algebra}, 86:187--221, 1993.

 \bibitem{Ma:introm}
S.~Majid.
\newblock Algebras and {H}opf algebras in braided categories.
\newblock volume 158 of {\em Lec. Notes in Pure and
Appl. Math}, pages 55--105. Marcel Dekker, 1994.

\bibitem{Ma:book}
S.~Majid.
\newblock {\em Foundations of Quantum Group Theory}.
\newblock Cambridge Univeristy Press, 1995.

\bibitem{JoyStr:bra}
A.~Joyal and R.~Street.
\newblock Braided monoidal categories.
\newblock Mathematics Reports 86008, Macquarie University, 1986.

\bibitem{Jac:int}
F.H. Jackson.
\newblock {$q$}-{I}ntegration.
\newblock {\em Proc. Durham Phil. Soc.}, 7:182--189, 1927.

\bibitem{Ma:fre}
S.~Majid.
\newblock Free braided differential calculus, braided binomial theorem and the
  braided exponential map.
\newblock {\em J. Math. Phys.}, 34:4843--4856, 1993.

\bibitem{KemMa:alg}
A.~Kempf and S.~Majid.
\newblock Algebraic $q$-integration and {F}ourier theory on quantum and braided
  spaces.
\newblock {\em J. Math. Phys.}, 35:6802--6837, 1994.
\end{thebibliography}

\end{document}